\RequirePackage[l2tabu, orthodox]{nag}
\RequirePackage{snapshot}

\documentclass[10pt,twocolumn]{article}

\makeatletter
\if@twocolumn
  \usepackage[dvips,letterpaper,top=0.5in, bottom=0.5in, left=0.75in, right=0.5in,includefoot,heightrounded]{geometry}
\else
  \usepackage[dvips,letterpaper,margin=1in,includefoot,heightrounded]{geometry}
\fi

\usepackage{srcltx}

\usepackage{amsmath}
\usepackage{amssymb,amsfonts}
\usepackage{amstext}

\usepackage{abstract}

\usepackage{graphics}
\usepackage{epsfig}
\usepackage[usenames,dvipsnames]{color}
\usepackage[usenames,dvipsnames]{xcolor}
\usepackage{subfigure}

\usepackage{booktabs}

\usepackage{setspace}
\usepackage{flushend}
\usepackage{multicol}

\usepackage{cite}
\usepackage{url}
\usepackage[normalem]{ulem}

\usepackage{enumerate}

\usepackage{hyperref}

\usepackage{lipsum}

\usepackage[sc,tiny]{titlesec}



\title{%
Wavelets for Elliptical Waveguide Problems
}

\author{%
M. M. S. Lira%
\thanks{%
M. M. S. Lira
is with
Power Systems Digital Laboratory,
Federal University of Pernambuco, P.O. Box 7800 Recife, Pernambuco, Brazil.
E-mail: \protect\url{milde@ufpe.br}
}
\quad
H. M. de Oliveira%
\thanks{%
H. M. de Oliveira
is with the
Communications Research Group,
Federal University of Pernambuco, P.O. Box 7800 Recife, Pernambuco, Brazil.
\protect\url{hmo@ufpe.br}
}
\quad
R. J. Cintra%
\thanks{%
R.~J.~Cintra is with
the Signal Processing Group,
Departamento de Estat\'{\i}stica,
Universidade Federal de Pernambuco,
Recife, Brazil.
This work was done while working
with the
Communications Research Group
at the same university.
E-mail:
\protect\url{rjdsc@stat.ufpe.org}
}
\quad
R. M. Campello de Souza%
\thanks{
R. M. Campello de Souza
is with the
Communications Research Group,
Federal University of Pernambuco, P.O. Box 7800 Recife, Pernambuco, Brazil.
\protect\url{ricardo@ufpe.br}
}
}

\date{}

\newcommand{\myabstract}{%
New elliptic cylindrical wavelets are introduced, which exploit the relationship between analysing filters and
Floquet's solution of Mathieu differential equations. It is shown that the transfer function of both multiresolution filters
is related to the solution of a Mathieu equation of odd characteristic exponent. The number of notches of these analysing
filters can be easily designed. Wavelets derived by this method have potential application in the fields of optics, microwaves
and electromagnetism.
}

\newcommand{\mykeywords}{%
Wavelets, Elliptic cylinder, Waveguide, Mathieu wavelets
}

\begin{document}

\makeatletter
\if@twocolumn

\twocolumn[%
  \maketitle
  \begin{onecolabstract}
    \myabstract
  \end{onecolabstract}
  \begin{center}
    \small
    \textbf{Keywords}
    \\\medskip
    \mykeywords
  \end{center}
  \bigskip
]
\saythanks

\else

  \maketitle
  \begin{abstract}
    \myabstract
  \end{abstract}
  \begin{center}
    \small
    \textbf{Keywords}
    \\\medskip
    \mykeywords
  \end{center}
  \bigskip
  \onehalfspacing
\fi

\section{Introduction}

A family of linear second-order differential equation with periodic coefficients was introduced by Emile Mathieu when studying vibrations in an elliptic drum~\cite{Mathieu}. Mathieu's equation is related to the wave equation for the elliptic cylinder. This paper is concerned with the canonical form of the Mathieu Equation:
\begin{equation}
\label{equation-1}
\frac{\operatorname{d}^{2}\!y}{\operatorname{d}\!w^2}
+
\big[
a-2q\cos\left ( 2w \right )
\big]
y=0.
\end{equation}
for
$a \in \mathbb{R}$, $q \in \mathbb{C}$.
The solution of \eqref{equation-1} is the elliptic cylindrical harmonic, known as Mathieu functions. Mathieu wavelets
as well as Mathieu functions can be appealing in a lot of Physics issues~\cite{Ruby} including vibrations in an elliptic drum, diffraction, amplitude distortion, the inverted pendulum, the radio frequency quadrupole, stability of a floating body, alternating gradient focusing, vibration in a medium with modulated density, and even when examining molecular dynamics of charged particles in electromagnetic traps~\cite{Nasse, Nie}. They have also long been applied on a broad scope of waveguide problems involving elliptical geometry~\cite{Holland, Sun, Shaw, Caorsi_et_al, Hussein, Wang, Schneider, Ragheb, Caorsi}, including:
(i) analysis for weak guiding for step index elliptical core optical fibres, (ii) power transport of elliptical waveguides, (iii) evaluating radiated waves of elliptical horn antennas, (iv) elliptical annular microstrip antennas with arbitrary eccentricity, and (v) scattering by a coated strip. In general, the solutions of \eqref{equation-1} are not periodic. However, for a given $q$, periodic solutions exist for infinitely many special values (eigenvalues) of $a$. For many physical solutions $y$ must be periodic of period $\pi$ or $2\pi$.  It is also convenient to distinguish even and odd periodic solutions, which are termed Mathieu functions of first kind.

\section{Preliminaries and Background}

One of four simpler types can be considered: Periodic solution ($\pi$ or $2\pi$) symmetry (even or odd). For $q \neq 0$,
the only periodic solution $y$ corresponding to any characteristic value $a=a_{\nu}(q)$ or $a=b_{\nu}(q)$ has the following
notation:
\begin{subequations}
\begin{align}
\emph{Even periodic solution:}
\nonumber
\\
\operatorname{ce}_{\nu}(w,q) &=
\sum_{m}^{} A_{\nu,m} \cos (mw),
\quad \text{for $a=a_{\nu}(q)$},
\\
\emph{Odd periodic solution:}
\nonumber
\\
\operatorname{se}_{\nu}(w, q) &=
\sum_{m}^{} A_{\nu,m} \sin (mw),
\quad \text{for $a= b_{\nu}(q)$},
\end{align}
\end{subequations}
where the sums are taken over even (respectively odd) values of $m$ if the period of $y$ is $\pi$ (respectively $2\pi$).
Given $\nu$, we denote henceforth $A_{\nu,m}$ by $A_m$, for short. Elliptic cosine and elliptic sine functions are represented
by $ce$ and $se$, respectively.
Interesting relationships are found when ${q \to 0}, \nu \neq 0$~\cite{Abramowitz}:
\begin{align}
\lim_{q \to 0} \operatorname{ce}_{\nu}(w,q) &= \cos(\nu w),
\\
\lim_{q \to 0} \operatorname{se}_{\nu}(w,q) &= \sin(\nu w).
\end{align}

As it happens with trigonometric functions, Mathieu functions hold orthogonality properties:

\begin{align}
\int_{0}^{2\pi} \operatorname{ce}_{\nu}(w,q) \operatorname{ce}_{\mu}(w,q)\operatorname{d}\!w& =0,
\quad
\nu \neq \mu,
\\
\int_{0}^{2\pi} \operatorname{se}_{\nu}(w,q) \cdot \operatorname{se}_{\mu}(w,q)\operatorname{d}\!w&=0,
\quad
\nu \neq \mu,
\\
\int_{0}^{2\pi} \operatorname{ce}_{\nu}(w,q) \cdot \operatorname{se}_{\mu}(w,q)\operatorname{d}\!w&=0.
\end{align}

The second non-periodic solution~\cite{Abramowitz} corresponding to $\operatorname{ce}_{\nu}(w,q)$ is the Mathieu function $Z\operatorname{ce}_{\nu}(w,q)$.
Similarly, The second non-periodic solution corresponding to $\operatorname{se}_{\nu}(w,q)$ is the Mathieu function $Z\operatorname{se}_{\nu}(w,q)$.

One of the most powerful results of Mathieu's functions is the Floquet's Theorem~\cite{McLachlan, Floquet}.
It states that periodic solutions of \eqref{equation-1} for any pair $(a,q)$ can be expressed in the form:
\begin{align}
y(w)&=F_{\nu}(w)= e^{j \nu w} \cdot P(w), \quad \text{or}
\\
y(w)&=F_{\nu}(-w)=e^{-j \nu w} \cdot P (-w),
\end{align}
where $\nu$ is a constant depending on $a$ and $q$ and $P(\cdot)$ is $\pi$-periodic in $w$. The constant $\nu$ is called the characteristic exponent. If $\nu$ is an integer, then $F_{\nu}(w)$ and $F_{\nu}(-w)$ are linear dependent solutions. Furthermore,
$y(w+k\pi) = e^{j \nu k\pi} y(w)$ or $y(w+k\pi)=e^{-j \nu k\pi} y(w)$, for the solution $F_{\nu}(w)$ or $F_{\nu}(-w)$, respectively.
We assume that the pair $(a,q)$ is such that $|\operatorname{cosh}(j \nu \pi)|<1$ so that the solution $y(w)$ is bounded on the real axis~\cite{Gradshteyn}.

The general solution of Mathieu's equation ($q \in R$, $\nu$ non-integer) has the form
\begin{equation}
y(w)=c_{1}e^{j \nu w} \cdot P(w) + c_{2}e^{-j \nu w} \cdot P(-w),
\end{equation}
where $c_1$ and $c_2$ are arbitrary constants.

All bounded solutions -those of fractional as well as integral order- are described by an infinite series of harmonic oscillations whose amplitudes decrease with increasing frequency. In the wavelet framework we are basically concerned with even solutions of period $2\pi$. In such cases there exist recurrence relations among the coefficients~\cite{Abramowitz}:
\begin{align}
(a-1-q)A_1 - qA_3&=0
,
\\
(a-m^2)A_m - q(A_{m-2} + A_{m+2}) &=0
,
\end{align}
for $m \geq 3$, $m$ odd.

\begin{figure}%
\centering
\includegraphics[width=0.4\textwidth]{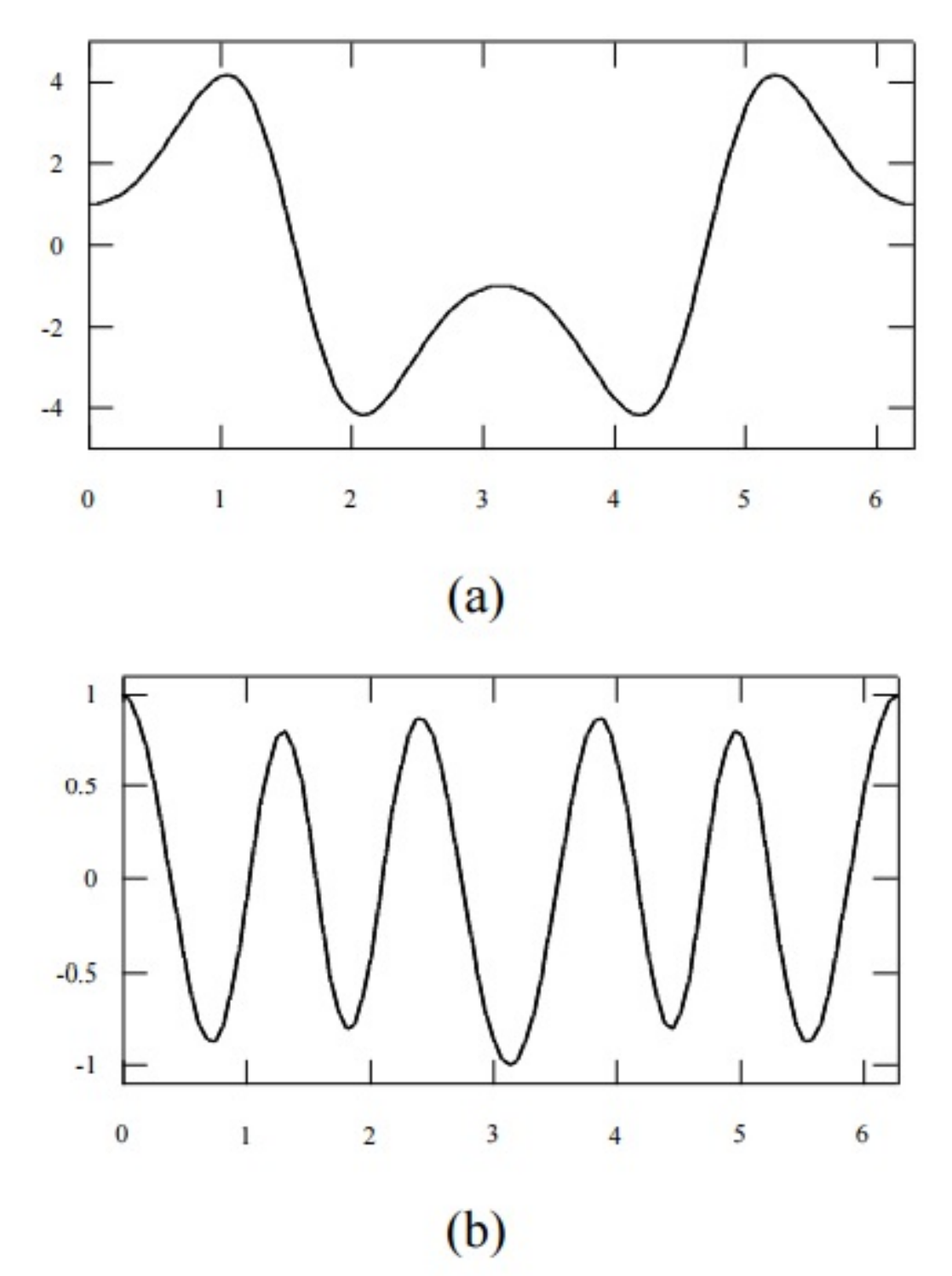}
\caption{ Some plots of $2\pi$-periodic $1^{\text{st}}$
kind even Mathieu functions.
Elliptic cosines shape for the following set
of parameters:
(a)~$\nu=1$ and $q=5$,
(b)~$\nu=5$ and $q=5$.}
\label{figure-1}
\end{figure}

Different methods for computing Mathieu functions are available in the literature~\cite{Frenkel,Alhargan}.
Figure~\ref{figure-1}
shows two illustrative waveforms of elliptic cosines, whose shape strongly depends on the parameters
$\nu$ and $q$.

\section{Elliptic Cylindrical Coordinates}

The $\nu$ coordinates are asymptotic angles of confocal hyperbolic cylinders, which are symmetric with respect to
the $x$ axis. The $u$ coordinates are confocal elliptic cylinders centred in the origin.
The transformation equations are given as follows:
\begin{align}
x &=a \operatorname{cosh}(u) \cos(v),
\\
y &=a \operatorname{sinh}(u) \sin(v),
\\
z &=z,
\end{align}
where $u\in [0, \infty)$, $v \in [0, 2\pi)$ and $z \in (-\infty, +\infty)$.

A further important property of Mathieu's functions is the orthogonality.
If $a(\nu+2\pi, q)$ and $a(\nu+2s, q)$ are
simple roots of $\cos\left ( \frac{\pi}{2}\nu \right )-y(\frac{\pi}{2})$, then~\cite{McLachlan}:
\begin{equation}
\int_{0}^{\pi} F_{\nu+2p}(w) \cdot F_{\nu+2s}(-w)\operatorname{d}\!w=0,
\quad
p \neq s.
\end{equation}
In other words, $\langle {F_{\nu+2p}(w),F_{\nu+2s}(-w)} \rangle  = 0$, $p \neq s$, where $\langle \cdot , \cdot \rangle$ denotes the inner product.

\section{Multiresolution Analysis Filters and Mathieu's Equation}

Wavelet analysis has quickly developed over the past years~\cite{Bulttheel} and there has been a explosion of wavelet applications including seismic geology, quantum physics, medicine, image processing (e.g., video data compression, reconstruction of high resolution images), fractals, computer graphics, linear system modelling, computer and human vision, denoising, filter banks, radar, wideband spreading, turbulence, statistics, volumetric visualisation, metallurgy, solution of partial differential equations, and so on. Essentially, the wavelet transform is a signal decomposition onto a set of basis functions, which is derived from a single prototype wavelet by scaling (dilatations and contractions) as well as translations (shifts).

In the sequel, wavelets are denoted by $\psi(t)$ and scaling functions by $\phi(t)$, with corresponding spectra $\Psi(w)$ and
$\Phi(w)$, respectively.
The equation $\sqrt{2}\sum_{n \in \mathbb{Z}} {h_n \phi(2t-n).}$, which is known as the \emph{dilation} or \emph{refinement equation}, is the chief relation determining a Multiresolution Analysis (MRA)~\cite{Mallat20, Mallat89}.

\section{Two Scale Relation of Scaling Function and Wavelet}

Defining the spectrum of the smoothing filter $\{h_k\}$ by
\begin{equation}
\label{equation-17}
H(w) \triangleq \frac{1}{\sqrt{2}} \sum_{k \in \mathbb{Z}} h_k e^{-jwk},
\end{equation}
the central equations (in the frequency domain)
of a multiresolution analysis are~\cite{Mallat89}:
\begin{align}
\Phi(w) = H \left (\frac{w}{2}  \right ) . \Phi \left (\frac{w}{2}  \right ),
\\
\Psi(w) = G\left (\frac{w}{2}  \right ) \Phi \left (\frac{w}{2}  \right ),
\end{align}
where
\begin{equation}
G(w) \triangleq \frac{1}{\sqrt{2}} \sum_{k \in \mathbb{Z}} g_k e^{-jwk},
\end{equation}
is the transfer function of the detail filter. The orthogonality condition corresponds to~\cite{Bulttheel, Mallat89}:
\begin{align}
H(0) = 1 \quad \text{and} \quad  H(\pi) = 0,
\\
|H(w)|^2 + |H(w+\pi)|^2=1,
\\
H(w) = -e^{-jw}G(w+\pi).
\end{align}

\section{Filters of the Mathieu MRA}

The subtle liaison between Mathieu's theory and wavelets was found by observing that the classical relationship~\cite{Bulttheel,Percival},
\begin{equation}
\Psi(w) = e^{-jw/2} H^{*}\left (\frac{w}{2}- \pi  \right ) \Phi\left (\frac{w}{2}  \right )
\end{equation}
presents a remarkable similarity to a Floquet's solution of a Mathieu's equation, since $H(w)$ is a periodic function.
As a first attempt, the relationship between the wavelet spectrum and the scaling function was put in the form:
\begin{equation}
\frac{ \Psi(w)} {\Phi\left (\frac{w}{2}\right ) }= e^{-jw/2}H^{*}\left (\frac{w}{2}- \pi  \right ).
\end{equation}

Here, on the second member, neither $\nu$ is an integer nor $H(\cdot)$ has a period $\pi$. By an appropriate scaling of
this equation, we can rewrite it according to:
\begin{equation}
\frac{\Psi(4w)}{\Phi(2w)}=e^{-j2w}.H\left ( 2w- \pi  \right ).
\end{equation}

It follows from \eqref{equation-17} that $H$ is periodic. We recognise that the function $Y(w) \triangleq\Psi(4w) / \Phi(2w)$ has
a nice interpretation in the wavelet framework. First, we recall that $\Psi(w) = G\left (\frac{w}{2} \right )\Phi \left (\frac{w}{2} \right )$
so that $\Psi(2w) =G(w) \Phi(w)$. Therefore the function related to Mathieu's equation is exactly $Y(w) = G(2w)$.
Introducing a new variable $z$, which is defined according to $2z \triangleq 2w-\pi$, it follows that $-Y\left ( z + \frac{\pi}{2}\right ) = e^{-j2z}H^{*}(2z)$.
The characteristic exponent can be adjusted to a particular value $\nu$
as follows:
\begin{equation}
-e^{-j( \nu-2)z}Y \left ( z+ \frac{\pi}{2} \right ) = e^{-j \nu z}.H(2z).
\end{equation}

Defining now $P(-z) \triangleq H^{*}(2z) = \sum_{k \in \mathbb{Z}} {c_{2k}e^{-j2kz}}$,
where $c_{2k} \triangleq  \frac{1}{\sqrt{2}} h^{*}_k$, we figure out that the right-side of the above equation represents a Floquet's solution
of some Mathieu differential equation. The function $P(\cdot)$ is $\pi$-periodic verifying the initial condition $P(0) =
 \frac{1}{\sqrt{2}} \sum_{k} {h_k}=1$, as expected. The filter coefficients are all assumed to be real.
Therefore, there exist a set of parameters $(a_G, q_G)$ such that the auxiliary function
\begin{equation}
y_{\nu}(z) \triangleq-e^{-j(\nu-2)z}Y_v \left ( z + \frac{\pi}{2}\right )
\end{equation}
is a solution of the following Mathieu equation:
\begin{equation}
\label{equation-29}
\frac{d^{2}y_{\nu}}{dz^{2}}+(a_G-2q_G \cos(2z))y_v = 0,
\end{equation}
subject to $y_{\nu}(0)=-Y(\pi/2) =-G(\pi) =-1$ and $\cos(\nu\pi)-y_{\nu}(\pi) = 0$, that is, $y_{\nu}(\pi) = (-1)^{\nu}$.

In order to investigate a suitable solution of \eqref{equation-29}, boundary conditions are established for predetermined $a,q$.
It turns out that when $\nu$ is zero or an integer, a belongs to the set of characteristic values $a_{\nu}(q)$.
The even ($2\pi$-periodic) solution of such an equation is given by:
\begin{equation}
y_{\nu}(z)=- \frac{\operatorname{ce}_{\nu}(z,q)}{\operatorname{ce}_{\nu}(0, q)}.
\end{equation}

The $Y_{\nu}(w)$ function associated with $y_{\nu}(z)$ and related to the detail filter of a `Mathieu MRA' is thus:
\begin{equation}
Y_{\nu}(w)=G_{\nu}(2w)= -e^{j(\nu-2)(w- \frac{\pi} {2})}.\frac{\operatorname{ce}_{\nu}(w- \frac{\pi} {2},q)} {\operatorname{ce}_{\nu}(0, q)}.
\end{equation}

Finally, the transfer function of the detail filter of a Mathieu wavelet is
\begin{equation}
G_{\nu}(w)= -e^{j(\nu-2)( \frac{w-\pi} {2})}.\frac{\operatorname{ce}_{\nu}( \frac{w-\pi} {2},q)} {\operatorname{ce}_{\nu}(0, q)}.
\end{equation}

The characteristic exponent $\nu$ should be chosen so as to guarantee suitable initial conditions, i.e., $G_{\nu}(0) = 0$
and $G_{\nu}(\pi)=1$, which are compatible with wavelet filter requirements~\cite{Bulttheel, Percival}. Therefore, $\nu$ must be odd.
It is interesting to remark that the magnitude of the above transfer function corresponds exactly to the modulus of an elliptic
sine~\cite{Gradshteyn}:
\begin{equation}
|G_{\nu}(w)| =\left |  \operatorname{se}_{\nu}
\left( \frac{w} {2},-q \right) / \operatorname{ce}_{\nu}(0, q) \right |.
\end{equation}

The solution for the smoothing filter $H(\cdot)$ can be found out via QMF conditions~\cite{Mallat20}, yielding:
\begin{equation}
H_{\nu}(w)= -e^{j\nu \frac{w} {2}}.\frac{\operatorname{ce}_{\nu}( \frac{w} {2},q)} {\operatorname{ce}_{\nu}(0, q)}.
\end{equation}
In this case, we find $H_{\nu}(\pi)=0$ and
\begin{equation}
|H_{\nu}(w)| = \left | \frac{\operatorname{ce}_{\nu}( \frac{w} {2},q)} {\operatorname{ce}_{\nu}(0, q)} \right |.
\end{equation}

Given $q$, the even first-kind Mathieu function with characteristic exponent $\nu$ is given by:
\begin{equation}
\operatorname{ce}_{\nu}(w,q) =\sum_{l=0}^{\infty} A_{2l+1} \cos(2l + 1)w,
\end{equation}
in which $\operatorname{ce}_{\nu}(0,q) =\sum_{l=0}^{\infty} A_{2l+1} $. The $G$ and $H$ filter coefficients of a Mathieu MRA can be expressed in terms
of the values $\{A_{2l+1}\}_{l \in \mathbb{Z}}$ of the Mathieu function as:
\begin{align}
\frac{h_{l}^{\nu}}{\sqrt{2}}&=- \frac{A_{|2l-\nu|}/2}{\operatorname{ce}_\nu(0, q)},
\\
\frac{g_{l}^{\nu}}{\sqrt{2}}&=(-1)^{l} \frac{A_{|2l+\nu-2|}/2}{\operatorname{ce}_\nu(0, q)}
.
\end{align}

\begin{figure}[!htb]
\centering
\includegraphics[width=0.4\textwidth]{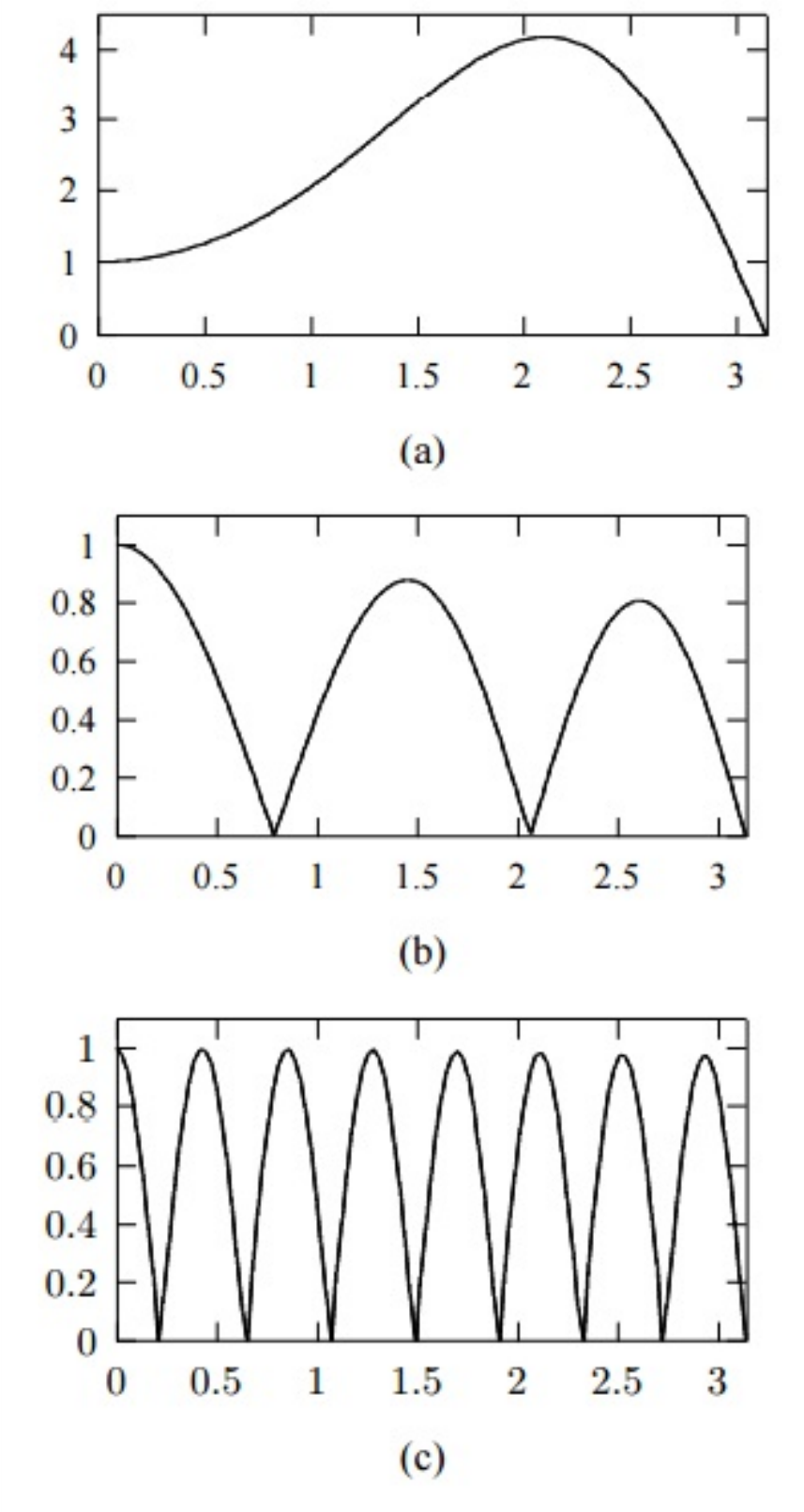}
\caption{Magnitude of the transfer function for Mathieu multiresolution analysis: smoothing filter $|Hv(w)|$ for a
few Mathieu parameters. (a) $\nu=1$, $q=5$, $a=1.858$; (b) $\nu=5$, $q=5$ $a=25.550$; (c) $\nu=15$,
$q =5$, $a=225.056$.}
\label{figure-3}
\end{figure}

It is straightforward to show that $h_{-l}^{\nu} = h_{l+\nu}^{\nu}$, and $\forall l \geq 0$. The normalising conditions are $\frac{1}{\sqrt{2}}  \sum_{k=-\infty }^{+\infty } h_{k}^{\nu} = -1$ and $\sum_{k=-\infty }^{+\infty } (-1)^{k}h_{k}^{\nu} = 0$.
Illustrative examples of filter transfer functions for a Mathieu MRA are shown in Figure~\ref{figure-3},
for $\nu=1,5$, and a particular value of $q$ (numerical solution obtained by 5-order Runge-Kutta method). The value of a is adjusted to an eigenvalue in each case, leading to a periodic solution.

Such solutions present a number of $\nu$ zeroes in the interval $|w|<\pi$. We observe lowpass behaviour (for the filter $H$)
and highpass behaviour (for the filter $G$), as expected.
Mathieu wavelets can be derived from the lowpass reconstruction filter by an iterative procedure (the same approach as the one usually used for plotting Daubechies wavelets). Infinite Impulse Response filters should be applied since Mathieu wavelet has no compact support.  However a Finite Impulse Response approximation can be generated by discarding negligible filter coefficients, say less than $10^{-10}$.

\begin{figure}[!ht]
\centering
\includegraphics[width=0.45\textwidth]{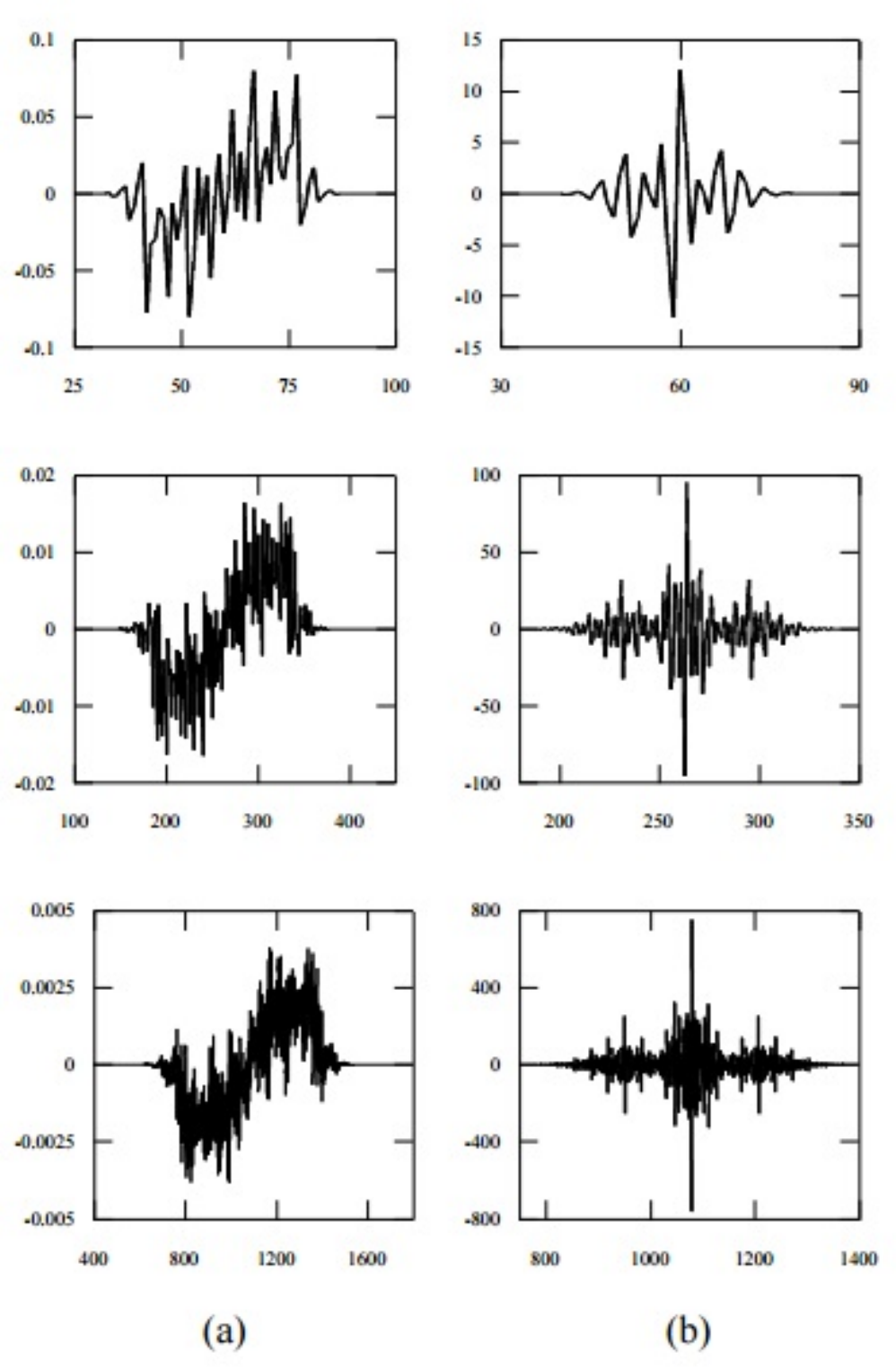}
\caption{FIR-Based Approximation of Mathieu Wavelets as the number of iteration increases: 2, 4, and 6 iterations, from the top to the bottom row.
Filter coefficients holding $|h| < 10^{-10}$ were thrown away (20 retained coefficients per filter in both cases).
(a)~Mathieu Wavelet with $\nu=1$ and $q=5$
and
(b)~Mathieu Wavelet with $\nu=5$ and $q=5$.}
\label{figure-4}
\end{figure}

In Figure~\ref{figure-4}, emerging pattern that progressively looks like the wavelet shape is shown. Waveforms were derived using the MATLAB wavelet toolbox. A fractal behaviour, which is common for some wavelets, can be noticed in these figures. As with many wavelets there is no nice analytical formula for describing Mathieu wavelets. Depending on the parameters $a$ and $q$ some waveforms (e.g., Figure.~\ref{figure-4}a) can present a somewhat unusual shape.

\section{Conclusions}

A new and wide family of elliptic cylindrical wavelets was introduced. It was shown that the transfer functions
of the corresponding multiresolution filters are related to Mathieu equation solutions. The magnitude of the detail
and smoothing filters corresponds to first-kind Mathieu functions with odd characteristic exponent. The number
of zeroes of the highpass $|G(w)|$ and lowpass $|H(w)|$ filters within the interval $|w|<\pi$ can be appropriately
designed by choosing the characteristic exponent. This seems to be the first connection found between Mathieu equations and wavelet theory. It opens new perspectives on linking wavelets and solutions of other differential equations (e.g. Associated Legendre functions, Coulomb wave function, Parabolic cylindrical functions etc.) Further generalisations such as the cases $q < 0$ or complex characteristic exponent can provide new interesting 'waves'. This new family of wavelets can particularly be an interesting tool for analysing optical fibres
due to its symmetry~\cite{Shaw, Wang, Vega}. Mathieu wavelets could as well be beneficial when examining molecular dynamics of charged particles in electromagnetic traps such as Paul trap or the mirror trap for neutral particles~\cite{Nasse, Nie}.

\section{Acknowledgments}

The authors thank Dr. Harold V. McIntosh from Departamento de Aplicacion de Microcomputadoras, Instituto de Ciencias,
Universidad Autonoma de Puebla, Mexico for information about Mathieu's functions.

\renewcommand\refname{References}

\end{document}